\documentclass[11pt]{article}
\usepackage{amssymb,latexsym}
\usepackage{srcltx}
\usepackage{graphicx}
\usepackage{amsmath}
\usepackage{amsfonts}
\usepackage{amscd}
\usepackage{mathrsfs}

\newcommand{\R}{\mathbb{R}}        
\newcommand{\N}{\mathbb{N}}        
\newcommand{\F}[1][]{Fr\'echet space#1}            
\newcommand{\FH}[1][]{Fr\'echet-Hilbert space#1}   

\newcommand{\Om}{{\rm (}\Omega{\rm )}}

\newcommand{\Lama}{\Lambda_\infty(\alpha)}

\newcommand{\qed}{\hspace*{\fill} $\Box $}

\newtheorem{theorem}{Theorem}[section]

\newtheorem{lemma}[theorem]{Lemma}

\def\proof{{\sc Proof.}\quad }

\begin{document}

\title{{\sc A characteristic property of the space $s$}}
\author{{\sc dietmar vogt}}

\date{}
\maketitle

\begin{abstract} It is shown that under certain stability conditions a complemented subspace of the space $s$ of rapidly decreasing sequences is isomorphic to $s$ and this condition characterizes $s$. This result is used to show that for the classical Cantor set $X$ the space $C_\infty(X)$ of restrictions to $X$ of $C^\infty$-functions on $\R$ is isomorphic to $s$, so completing the theory developed in \cite{VP}.
\end{abstract}

\footnotetext{\hskip -.8em   {\em 2010 Mathematics Subject Classification.}
    {46A45, 46A63, 46E10}\hfil\break\indent
{\em Key words and phrases: space $s$, stability condition, Cantor set}  \hfil\break\indent
    }

\section{Introduction}

In the present note we study the space $s$ of rapidly decreasing sequences, that is, the space
$$s=\{x=(x_0,x_1,\dots)\,:\, \lim_n x_n n^k=0 \text{ for all }k\in\N\}.$$
Equipped with the norms $\|x\|_k=\sup_n|x_n|(n+1)^k$ it is a nuclear \F. It is isomorphic to many of the Fr\'echet spaces which occur in analysis, in particular, spaces of $C^\infty$-functions.

It is easily seen that instead of the sup-norms we might use the norms
$$|x|_k=\big(\sum_n |x_n|^2(n+1)^{2k}\big)^{1/2}$$
which makes $s$ a \FH.

More generally, we define for any sequence $\alpha: 0\le\alpha_0\le\alpha_1\le\nearrow +\infty$ the power series space of infinite type
$$\Lama:=\{x=(x_0,x_1,\dots)\,:\,|x|_t^2=\sum_{n=0}^\infty |x_n|^2 e^{2t\alpha_n}<\infty\text{ for all }t>0\}.$$
Equipped with the hilbertian norms $|\cdot|_k$, $k\in\N_0$, it is a \FH. It is nuclear if, and only if, $\limsup_n\log n/\alpha_n<\infty$. With this definition $s=\Lama$ with $\alpha_n = \log (n+1)$.

A \F{} with the fundamental system of seminorms $\|\cdot\|_0\le\|\cdot\|_1\le\dots$ has property (DN) if
$$\exists p\, \forall k\, \exists K,\,C>0: \|\cdot\|_k^2\le C \|\cdot\|_p\; \|\cdot\|_K.$$
In this case $\|\cdot\|_p$ is called a dominating norm.

$E$ has property $\Om$ if
$$\forall p\, \exists q\, \forall m\,\exists 0<\theta<1, C>0: \|\cdot\|_q^*\le C {\|\cdot\|_p^*}^\theta\; {\|\cdot\|_m^*}^{1-\theta}.$$
Here we set for any continuous seminorm $\|\cdot\|$ and $y\in E'$ the dual, extended real valued, norm $\|y\|^*=\sup\{|y(x)|\,:\,x\in E,\;\|x\|\le 1\}$.

By  Vogt-Wagner \cite{VW1}  a \F{} $E$ is isomorphic to a complemented subspace of $s$ if, and only if, it is nuclear and had properties (DN) and $\Om$.

It is a long standing unsolved problem of the structure theory of nuclear \F s, going back to Mityagin, whether every complemented subspace of $s$ has a basis. If it has a basis then it is isomorphic to some power series space $\Lama$. The space $\Lama$ to which it is isomorphic, if it has a basis, can be calculated in advance by a method going back to Terzio\u{g}lu \cite{Terz1} which we describe now.

Let $X$ be a vector space and $A\subset B$ absolutely convex subsets of $X$. We set
$$\delta_n(A,B):=\inf\{\delta>0\,:\,\text{exists linear subspace }F\subset X, \dim F\le n\text{ with }A\subset \delta B+ F\}.$$
It is called the $n$-th Kolmogoroff diameter of $A$ with respect to $B$.

If now $E$ is a complemented subspace of $s$, that is, $E$ is nuclear and has properties (DN) and $\Om$, then we choose $p$ such that $\|\cdot\|_p$ is a dominating norm and for $p$ we choose $q>p$ according to property $\Om$. We set
$$\alpha_n=-\log \delta_n(U_q,U_p)$$
where $U_k=\{x\in E\,:\,\|x\|_k\le1\}.$ The space $\Lama$ is called the associated power series space and $E\cong \Lama$ if it has a basis.

If $\limsup_n \alpha_{2n}/\alpha_n<\infty$ then, by  Aytuna-Krone-Terzio\u{g}lu  \cite[Theorem 2.2]{AKT}, $E\cong\Lama$. This is, in particular, the case if $E$ is stable, that is, if $E\oplus E\cong E$.

For all that and further results of structure theory of infinite type power series spaces see \cite{V7}, for results and unexplained notation of general functional analysis see \cite{MV}.

\section{Main result}

\begin{lemma}\label{lem1} Let $E$ be a complemented subspace of $s$, $\|\cdot\|_0$ a dominating hilbertian norm and $\|\cdot\|_1$ a hilbertian norm chosen for $\|\cdot\|_0$ according to $\Om$. If there is a linear isomorphism $\psi: E\oplus E\to E$ such that
 \begin{eqnarray*}\|x\|_0+\|y\|_0&\le &C_0\, \|\psi(x\oplus y)\|_0\\ \|\psi(x\oplus y)\|_1&\le& C_1 (\|x\|_1+\|y\|_1)\end{eqnarray*}
 then $E\cong s$.
 \end{lemma}

 \proof  For $x\oplus y\in E\oplus E$ we set $|||(x,y)|||_0:= (\|x\|^2_0+\|y\|^2_0)^{1/2}$ and $|||(x,y)|||_1:= (\|x\|^2_1+\|y\|^2_1)^{1/2}$. With new constants $C_k$ we have
 \begin{equation}\label{eq1}|||x\oplus y|||_0\le C_0 \|\psi(x\oplus y)\|_0\text{ and }\|\psi(x\oplus y)\|_1\le C_1 |||x\oplus y|||_1.\end{equation}

 To calculate the associated power series space for $E$ we set:
 \begin{eqnarray*}
 \alpha_n &=& -\log\delta_n(U_1,U_0)\text{ where }U_k=\{x\in E\,:\,\|x\|_k\le 1\},\\
 \beta_n &=& -\log\delta_n(V_1,V_0)\text{ where }V_k=\{x\oplus y\in E\oplus E\,:\,|||x\oplus y|||_k\le 1\}.
 \end{eqnarray*}
 Due to the estimates (\ref{eq1}) we have
 $$\frac{1}{C_1}\psi(V_1)\subset U_1\subset U_0\subset C_0\psi(V_0)$$
  and therefore
 $$\delta_n(V_1,V_0)= \delta_n(\psi V_1,\psi V_0) \le C_0 C_1 \delta_n(U_1,U_0)$$
which implies
 $$\alpha_n \le \beta_n+d$$
 with $d=\log C_0 C_1$.

 By explicit calculation of the Schmidt expansion of the canonical map $j_1^0$ between the local Hilbert spaces of $|||\cdot|||_1$ and $|||\cdot|||_0$ and by use of the fact that singular numbers and Kolmogoroff diameters coincide, we obtain
 that $\beta_{2n}=\beta_{2n+1}=\alpha_n$ for all $n\in \N_0$.

 Therefore we have $\alpha_{2n}´\le\beta_{2n}+d= \alpha_n+d$ for all $n\in\N_0$ and this implies $\alpha_{2^k}\le \alpha_1 + k\,d$ for all $k\in\N_0$. For $n\in\N$ we find $k\in\N$ such that $2^{k-1}\le n \le2^k$ and we obtain $\alpha_n\le \alpha_{2^k}\le \alpha_1 + k\,d\le (\alpha_1+d) +d\log n$.

 Since $E\subset s$, which implies the left inequality below, we have shown that there is a constant $D>0$ such that
 $$\frac{1}{D}\log n\le \alpha_n\le D \log n$$
 for large $n\in\N$. This implies that $\Lambda_{\infty}(\alpha)=s$. \qed

 A \FH{} $E$ is called \em normwise stable \rm if it admits a fundamental system of hilbertian seminorms for which there is an isomorphism $\psi: E\oplus E\to E$ such that
 $$\frac{1}{C_k}(\|x\|_k+\|y\|_k)\le \|\psi(x\oplus y)\|_k\le C_k (\|x\|_k+\|y\|_k)$$ for all $k$. Since, clearly, $s$ is normwise stable we have shown.
 \begin{theorem} $E\cong s$ if, and only if, $E$ is isomorphic to a complemented subspace of $s$ and normwise stable. \end{theorem}
 We may express Lemma \ref{lem1} also in the following way:
 \begin{theorem}\label{thm2} Let the \FH{} $E$ be a complemented subspace of $s$, $\|\cdot\|_0$ a dominating norm and $\|\cdot\|_1$ be a norm chosen according to $\Om$. Let $P$ be a linear projection in $E$, continuous with respect to $\|\cdot\|_0$. We set $E_1=R(P)$, $E_2=N(P)$ and assume that there are linear isomorphisms $\psi_j:E\to E_j$, $j=1,2$, continuous with respect to $\|\cdot\|_1$ such that $\psi^{-1}$ is continuous with respect to $\|\cdot\|_0$. Then $E\cong s$.
 \end{theorem}

 \proof We set $\psi(x\oplus y):=\psi_1(x)+\psi_2(y)$ and obtain with suitable constants:
 \begin{eqnarray*} \|x\|_0+\|y\|_0 &\le& C'(\|\psi_1(x)\|_0+\|\psi_2(y)\|_0)
 \le C_0\|\psi_1(x)+\psi_2(y)\|_0=C_0\|\psi(x\oplus y)\|_0\\
 \|\psi(x\oplus y)\|_1&=&\|\psi_1(x)+\psi_2(y)\|_1\le\|\psi_1(x)\|_1+\|\psi_2(y)\|_1\le C_0(\|x\|_1+\|y\|_1).
 \end{eqnarray*}
 Lemma \ref{lem1} yields the result. \qed

 \section{Application}

 An interesting application of this result is the following. Let $X\subset [0,1]$ be the classical Cantor set and $C_\infty(X):=\{f|_X\,:\, f\in C^\infty[0,1]\}=\{f|_E\,:\,f\in C^\infty(\R)\}.$ The space $C_\infty(X)$ equipped with the quotient topology is a nuclear \F{} and, since $C^\infty[0,1]\cong s$ isomorphic to a quotient of $s$, hence has property $\Om$. By a theorem of Tidten \cite{T2} it has also property (DN). Therefore it is isomorphic to a complemented subspace of $s$ (see \cite{VW1}).

 We should remark that, due to the fact that $X$ is perfect, we have $C_\infty(X)=\mathcal{E}(X)$ where $\mathcal{E}(X)$ denotes the space of Whitney jets on $X$, for which Tidten's result is formulated.

 By obvious identifications we have
 $$C_\infty(X)\cong C_\infty(X\cap[0,1/3])\oplus C_\infty(X\cap[2/3,1])\cong C_\infty(X)\oplus C_\infty(X)$$
 and it is easily seen that this establishes normwise stability. Therefore we have shown
 \begin{theorem} If $X$ is the classical Cantor set, then $C_\infty(X)\cong s$. \end{theorem}

 It should be remarked that in \cite{AGK} it has been shown that for the Cantor set $X$ the diametral dimensions of $\mathcal{E}(X)$ and $s$ coincide, from where, by means of the Aytuna-Krone-Terzio\u glu Theorem, on can derive the same result.

 Referring to the terminology of \cite{VP} we have also shown that $A_\infty(X)\cong s$ which completes the theory developed in \cite{VP}.


\vspace{.5cm}

\noindent Bergische Universit\"{a}t Wuppertal,
\newline FB Math.-Nat., Gau\ss -Str. 20,
\newline D-42119 Wuppertal, Germany
\newline e-mail: dvogt@math.uni-wuppertal.de

\end{document}